\documentstyle[12pt]{article}
\title{A simple proof of exponential decay in the  two dimensional percolation model
} 
\author{Yu Zhang\footnote{Research supported by NSF grant DMS-0405150. }}
\begin{document}
\baselineskip .20in
\maketitle
\begin{abstract}
In 1980,
Kesten showed the exponential decay of percolation probability in the subcritical phase
for the two-dimensional percolation model.   
This result implies his celebrated computation
that $p_c=0.5$ for bond percolation in the square lattice, and site percolation in the triangular lattice,
respectively.
In this paper, we present a simpler proof for Kesten's theorem.
\end{abstract}

\section{ Introduction and statement of results.}
We may deal with the percolation model on the two-dimensional periodic lattice (see the detailed definition in Kesten (1982)).
For simplicity, we select the  triangular lattice without loss of generality,
since we do not need to deal with the dual lattice separately. 
Consider site percolation on the triangular lattice. Each vertex of the
lattice is open with probability $p$ and closed with probability $1-p$, and the sites are open
independently of each other. We will realize the triangular lattice with vertex set $
\mathbf{Z}^{2}$. For a given $(x,y)\in \mathbf{Z}^{2}$, its nearest
neighbors are defined as $(x\pm 1,y)$, $(x,y\pm 1)$, $(x+1,y-1)$, and $
(x-1,y+1)$. Bonds between neighboring or adjacent sites therefore correspond to
vertical or horizontal displacements of one unit, or diagonal displacements
between the two nearest vertices along a line making an angle of $135^{\circ }$
with the positive $x$-axis. 
Recall that the triangular lattice may also be viewed with sites as hexagons
in a regular hexagonal tiling of the plane. 
The corresponding probability measure on the {\em configurations} of open 
and closed  sites  is denoted by ${\bf P}_{p}$. 
We also denote by ${\bf E}_{p}$ the expectation with respect to ${\bf P}_{p}$.

A path  from $u$ to $v$ is a sequence $(v_0,e_1, v_1,...,v_{i}, e_{i+1}, v_{i+1},...,v_l)$
with distinct vertices $v_i$ ($0\leq i\leq n$) and $v_0=u$ and $v_n=v$ and with bonds $e_i$  between $v_i$ and $v_{i+1}$.
If  $u=v$, the path is called a {\em circuit}.
If all of the sites in a path are open, the path is called an {\em open path}.
Given a rectangle $[-n,n]\times [-m, m]$, a {\em left-right open  crossing} is path $(v_0,e_1, v_1,...,v_{i}, e_{i+1}, v_{i+1},...,v_l)$ in $[-n,n]\times [-m, m]$ such that all of its vertices  inside $(-n,n)\times (-m,m)$
are open except $v_0$ and $v_l$,
which are at the left and at the right of the rectangle, respectively. 
Similarly, we can define a {\em top-bottom open
crossing}. We denote by $LR([-n,n]\times [-m, m])$ and $TB([-n,n]\times [-m, m])$ the events that
there exist the left-right and the top-bottom open crossings, respectively. We may replace
the open paths with  closed paths to have the events of $LR^*([-n,n]\times [-m, m])$ and $TB^*([-n,n]\times [-m, m])$,
respectively. We call them the {\em left-right closed crossing} and the {\em top-bottom closed crossing}, respectively.

There might be many  open crossings. For each open crossing
$\Gamma$,  it divides $[-n,n]\times [-m, m]$ into two parts:
the top part $T(\Gamma)$ and the bottom part $B(\Gamma)$, including the top and the bottom of $[-n,n]\times [-m,m]$,
respectively. 
We also denote by
$T^o(\Gamma)$ and $B^o(\Gamma)$ the interiors of 
$T(\Gamma)$ and $B(\Gamma)$, respectively. If there is more than one  left-right open crossing,
 we select an open crossing with the smallest vertex set
$B(\Gamma)$. We call the open crossing  the {\em lowest crossing}. Without loss of generality, we still
denote  by $\Gamma$ the lowest open crossing.
By the definition of the lowest crossing, 
it can be obtained (see page 317 in Grimmett (1999) or Proposition 2.3 in Kesten (1982))
 that the event of $\{\Gamma=\gamma\}$ for some fixed
left-right crossing $\gamma$ only depends on the open or closed vertices  in $B(\gamma)$.
This property is said to be  the {\em independent property} of the lowest crossing.
%Similarly, we can define the most left top-bottom crossing. 

If $\Gamma$ is the lowest open crossing, then for each $v\in \Gamma$, 
it is well known (see Proposition 2.2 in Kesten (1982))
that there exists a closed path from $v$ (not including  $v$)
to the bottom of $[-n,n]\times [-m,m]$. 
%The proof of this argument can be found in the Appendix of Kesten
%(1982). Let us describe the proof briefly. 
%For each $v\in \Gamma$, let
%${\bf C}^*(v)$ be the closed cluster that is adjacent to $v$ inside $B(\Gamma)$. 
%By Proposition 2.2 in Kesten (1982), the boundary of the cluster is a circuit. The circuit may contains
%vertices of $\Gamma$ and the bottom of $[-n,n]\times [-m,m]$. The other vertices of the circuit 
%except the bottom are open.
%If ${\bf C}^*(v)$
%does not contain a vertex of the bottom of $[-n,n]\times [-m,m]$, then we can use the boundary of 
%${\bf C}^*(v)$ and $\Gamma$ to construct another open path which is ``lower" than $\Gamma$.
%This contradicts the assumption that $\Gamma$ is the lowest.
By this observation, for each $v\in \Gamma$, there exist two disjoint open paths 
$\Gamma_1$ and $\Gamma_3$  with 
$$\Gamma_1\cup \{v\}\cup \Gamma_3=\Gamma,$$
from $v$ to the left and to the right of $[-n,n]\times [-m,m]$, respectively. In addition, there exists a closed path
$\Gamma_2$ in $B^\circ(\Gamma)$ from $v$ to the bottom of $[-n,n]\times [-m,m]$. 
On the other hand, by using Proposition 2.2 in Kesten (1982) again, if there exist the three paths
$\Gamma_i$ for $i=1,2,3$ at $v$, and $v$ is open, then $v$ is on the lowest crossing.
This is called the 
{\em three-arm-path argument} for each vertex on the lowest crossing $\Gamma$. Similarly, the
three-arm-path argument holds for the 
left-most top-bottom crossing.

We may generalize these arguments into a circuit enclosed by a path. Let 
$\Delta$ be an open set surrounded by a circuit  $\partial \Delta$.  
We select four vertices $v_i$ for $i=1,2,3,4$ from the circuit. 
Let $L$ (left), $T$ (top), $R$ (right), and $B$ (bottom) be the paths along $\Delta$ clockwise 
from $v_1$ to $v_2$, from $v_2$ to $v_3$, from $v_3$ to $v_4$,
and from $v_4$ to $v_1$, respectively.  With these paths,  we can define the events $LR(\Delta)$ and
$TB^*(\Delta)$ in the same way that we did for a rectangle. In addition, we can also define the lowest 
left-right,  and the  left-most  top-bottom open or closed crossings. We denote by $\Gamma_{LR}(\Delta)$ and $\Gamma_{TB}^*(\Delta)$
the lowest left-right open, and the  left-most  top-bottom closed crossings. By the same discussion above, 
the independent property of the lowest crossing and the three-arm-path argument still hold for $\Gamma_{LR}(\Delta)$ and $\Gamma_{TB}^*(\Delta)$.
In fact, Kesten (1982) discussed these topology properties precisely in the circuit as we defined above.

For each $v\in [-n,n]^2$, we say there are three arm paths from $v$, as we mentioned,
if there are two disjoint  open paths $\Gamma_1$ and $\Gamma_3$ in $[-n,n]^2$
from $v$ to the left and to the right of $[-n,n]^2$, and there exists a closed path $\Gamma_2$ from $v$ to the bottom
of $[-n,n]^2$, respectively.  Also, $v$ is open.
Moreover, if there exists an additional  closed path $\Gamma_4$ in $T^\circ (\Gamma)$ from 
$v$ to the top boundary of $[-n,n,]^2$, 
 we say
there are {\em four arm paths}  at $v$.
If there are four arm paths at $v$, $v$ is said to be a {\em pivotal vertex} of open crossing $LR([-n,n]^2)$.
Let $N_n$ be  all the  pivotal vertices in $[-n,n]^2$. 

 The {\em open cluster} of the vertex $x$, 
${\bf C}(x)$,
consists of all vertices, that including $x$, that are connected  by  open paths.
For any collection ${\bf A}$ of vertices, $|{\bf A}|$  
denotes the cardinality of $A$. We choose ${\bf 0}$ as the origin. The percolation probability and the mean size
of the open cluster are denoted by
\begin{eqnarray*}
\theta (p)= {\bf P}_{p}(|{\bf C}({\bf 0})|=\infty)\mbox{ and } \chi(p)= {\bf E}_p (|{\bf C}({\bf 0})|),
\end{eqnarray*}
and the critical probabilities are defined by
$$p_c=\sup\{p:\theta (p)=0\}\mbox{ and }p_T= \sup\{p: \chi(p) < \infty\}.$$
Similarly, we denote by ${\bf C}^*(x)$ the {\em closed cluster} including $x$.
With these definitions, the crucial step in Kesten's paper
(1980) is to estimate ${\bf E}_{p} (N_n, LR([-n,n]^2))$ in the following theorem. 
In this paper, we will present a simpler proof for his estimate.\\

{\bf Theorem 1.} {\em If $p \leq 0.5$, then there exists $\alpha >0$ such that for all $n$,}
$$n^\alpha\leq {\bf E}_{p}\left(N_n\,\,|\,\, LR\left([-n,n]^2\right) \right) .\eqno{(1.1)}$$

With Theorem 1, we will have the following corollary.\\

{\bf Corollary 2.} {\em If $p < 0.5$, then there exist constants $C_i=C_i(p)$ for $i=1,2$  such that
$${\bf P}_p ({\bf C}({\bf 0})\cap \partial [-n,n]^2 \neq \emptyset )\leq C_1 \exp(-C_2 n)\eqno{(1.2)}$$
and }
$${\bf P}_{1-p} ({\bf C}^*({\bf 0})\cap \partial [-n,n]^2 \neq \emptyset )\leq C_1 \exp(-C_2 n),\eqno{(1.3)}$$
where $\partial [-n,n]^2$ is the boundary vertex set of $[-n,n]^2$.\\

{\bf Remark.} If (1.2) holds, it follows from Theorem 5.1 in Kesten (1982) that 
$${\bf P}_p (|{\bf C}({\bf 0})|\geq n )\leq C_1 \exp(-C_2 n).$$

For more than two decades since 1959, one of the most important discoveries in the history of the percolation model
was the rigorous determination of $p_c=0.5$ for the square lattice and the triangular lattice. 
 Harries (1960) proved that $p_c \geq 0.5$. The precise lower bound of $p_c$ seems much harder to achieve.
After 20 years, by the estimate in Theorem 1, Kesten (1980) finally showed  that $p_c=0.5$.
In this paper, we present a proof by using Theorem 1 to show $p_c=0.5$.\\

{\bf Corollary 3.} $p_c=p_T=0.5$.\\

{\bf Remark.} The same argument can be carried out to show Theorem 1 if $p \leq p_c$, and Corollary 2 if $p < p_c$
for the percolation model in the two-dimensional periodic lattice.
In addition, the same argument can also be carried out to show that $p_c=p_T=0.5$ for the bond percolation model
in the square lattice.\\

\section{ Proofs of theorems and corollaries.}
Before the proofs of the  Theorems, we introduce a lemma by Russo (1978) and Seymour and Welsh (1978).\\

{\bf RSW lemma.} {\em If ${\bf P}_p (LR([0,n]^2)) \geq \delta > 0$, then for each integer $k$,
there exists a positive constant $C_3=C_3(k, \delta)$ such that}
$${\bf P}_{p} \left ( LR([0,kn]\times [0, n])\right) \geq C_3.\eqno{(2.1)}$$

Note that by symmetry, we know that 
$${\bf P}_{0.5} \left(LR([0,n]^2)\right) =0.5\mbox{ and } {\bf P}_{0.5} \left(LR^*([0,n]^2)\right) =0.5$$
 for each $n$, so by the RSW lemma,
$${\bf P}_{0.5} ( LR([0,kn]\times [0, n])) \geq C_3 \mbox{ and } {\bf P}_{0.5} ( LR^*([0,kn]\times [0, n])) \geq C_3.\eqno{(2.2)}$$

With (2.2) and the FKG inequality,  we can directly show that there exists a closed circuit in an annulus
with a positive probability. Thus the following lemma 
can  be directly obtained by this probability estimate (see Theorem 11.89 in Grimmett (1999)). \\

{\bf Lemma 1.} {\em There exists $C_4>0$ such that} 
$${\bf P}_{0.5} \left( {\bf C}({\bf 0}) \cap \partial [-n, n]^2\neq \emptyset \right) \leq n^{-C_4}.$$

Now we show Theorem 1 by using Lemma 1.\\
\begin{figure}
\begin{center}
\setlength{\unitlength}{0.0125in}%
\begin{picture}(400,150)(107,740)
\thicklines
%\put(-60,770){\framebox(200,200)[br]{\mbox{$[-2n,2n]^2$}}}
\put(220,770){\framebox(200,200)[br]{\mbox{$$}}}
\put(270,820){\framebox(100,100)[br]{\mbox{$$}}}

\put(220,720){\line(0,1){200}}
\put(270,720){\line(0,1){250}}
\put(320,720){\line(0,1){250}}

\put(370,720){\line(0,1){250}}
\put(420,720){\line(0,1){250}}

\put(220,720){\line(1,0){200}}
\put(220,920){\line(1,0){200}}

\put(370,980){$n/2$}
\put(250,980){$-n/2$}
\put(320,980){$0$}
\put(370,770){\line(0,1){200}}

\put(220,890){\circle*{4}}
\put(225,890){\circle*{4}}
\put(230,890){\circle*{4}}
\put(235,890){\circle*{4}}
\put(240,895){\circle*{4}}
\put(245,900){\circle*{4}}
\put(250,900){\circle*{4}}
\put(255,900){\circle*{4}}

\put(260,900){\circle*{4}}
\put(265,900){\circle*{4}}
\put(270,900){\circle*{4}}
\put(275,900){\circle*{4}}
\put(280,900){\circle*{4}}
\put(285,900){\circle*{4}}
\put(290,900){\circle*{4}}
\put(295,900){\circle*{4}}
\put(300,900){\circle*{4}}
\put(305,900){\circle*{4}}
\put(310,900){\circle*{4}}
\put(315,900){\circle*{4}}
\put(315,895){\circle*{4}}
\put(315,890){\circle*{4}}
\put(315,885){\circle*{4}}
\put(320,885){\circle*{4}}
\put(325,885){\circle*{4}}
\put(330,885){\circle*{4}}
\put(335,885){\circle*{4}}
\put(340,885){\circle*{4}}
\put(345,885){\circle*{4}}
\put(350,885){\circle*{4}}
\put(355,885){\circle*{4}}
\put(360,885){\circle*{4}}
\put(365,885){\circle*{4}}
\put(370,885){\circle*{4}}
\put(375,885){\circle*{4}}
\put(380,885){\circle*{4}}
\put(385,885){\circle*{4}}
\put(390,885){\circle*{4}}
\put(395,885){\circle*{4}}
\put(400,885){\circle*{4}}
\put(405,885){\circle*{4}}
\put(410,885){\circle*{4}}
\put(415,885){\circle*{4}}
\put(420,885){\circle*{4}}
\put(200,910){${n}$}
\put(200,970){$3n\over 2$}
\put(200,950){$L$}
\put(430,950){$R$}
\put(300,975){$T$}

\put(220,950){\circle*{4}}
\put(225,950){\circle*{4}}
\put(230,950){\circle*{4}}
\put(235,950){\circle*{4}}
\put(240,950){\circle*{4}}
\put(245,950){\circle*{4}}
\put(250,950){\circle*{4}}
\put(255,950){\circle*{4}}
\put(260,950){\circle*{4}}
\put(265,950){\circle*{4}}
\put(270,950){\circle*{4}}
\put(275,950){\circle*{4}}
\put(280,950){\circle*{4}}
\put(285,950){\circle*{4}}
\put(270,950){\circle*{4}}
\put(275,950){\circle*{4}}
\put(280,950){\circle*{4}}
\put(285,950){\circle*{4}}
\put(270,950){\circle*{4}}
\put(275,950){\circle*{4}}
\put(280,950){\circle*{4}}
\put(285,950){\circle*{4}}
\put(290,950){\circle*{4}}
\put(295,950){\circle*{4}}
\put(300,950){\circle*{4}}
\put(305,950){\circle*{4}}
\put(310,955){\circle*{4}}
\put(315,955){\circle*{4}}
\put(320,955){\circle*{4}}

%\put(290,975){\circle*{4}}
\put(290,970){\circle*{4}}
\put(290,965){\circle*{4}}
\put(290,960){\circle*{4}}
\put(290,955){\circle*{4}}
\put(290,950){\circle*{4}}
\put(290,945){\circle*{4}}
\put(290,940){\circle*{4}}
\put(290,935){\circle*{4}}
\put(290,930){\circle*{4}}
\put(290,925){\circle*{4}}
\put(290,920){\circle*{4}}
\put(285,915){\circle*{4}}
\put(280,910){\circle*{4}}
\put(275,905){\circle*{4}}
\put(270,900){\circle*{4}}

%\put(350,975){\circle{3}}
\put(350,970){\circle{3}}
\put(350,965){\circle{3}}
\put(350,960){\circle{3}}
\put(350,955){\circle{3}}
\put(350,950){\circle{3}}
\put(350,945){\circle{3}}
\put(350,940){\circle{3}}
\put(350,935){\circle{3}}
\put(350,930){\circle{3}}
\put(350,925){\circle{3}}
\put(350,920){\circle{3}}
\put(350,915){\circle{3}}
\put(350,910){\circle{3}}
\put(350,905){\circle{3}}
\put(350,900){\circle{3}}
\put(350,895){\circle{3}}
\put(350,890){\circle{3}}
\put(350,885){\circle{3}}

\put(345,920){\circle{3}}
\put(340,915){\circle{3}}
\put(335,915){\circle{3}}
\put(330,915){\circle{3}}
\put(325,910){\circle{3}}
\put(320,910){\circle{3}}
\put(315,910){\circle{3}}
\put(310,910){\circle{3}}
\put(305,905){\circle{3}}
\put(305,900){\circle{3}}

\put(300,905){\circle*{4}}
\put(300,910){\circle*{4}}
\put(300,915){\circle*{4}}
\put(300,920){\circle*{4}}
\put(300,925){\circle*{4}}
\put(300,930){\circle*{4}}
\put(300,935){\circle*{4}}
\put(295,935){\circle*{4}}

\put(305,895){\circle{3}}
\put(305,890){\circle{3}}
\put(305,885){\circle{3}}
\put(305,880){\circle{3}}
\put(305,875){\circle{3}}
\put(305,870){\circle{3}}
\put(305,865){\circle{3}}
\put(305,860){\circle{3}}
\put(305,855){\circle{3}}
\put(305,850){\circle{3}}
\put(305,845){\circle{3}}
\put(305,840){\circle{3}}
\put(305,835){\circle{3}}
\put(305,830){\circle{3}}
\put(305,825){\circle{3}}
\put(305,820){\circle{3}}
\put(300,815){\circle{3}}
\put(295,810){\circle{3}}
\put(290,805){\circle{3}}
\put(290,800){\circle{3}}
\put(290,795){\circle{3}}
\put(290,790){\circle{3}}
\put(290,785){\circle{3}}
\put(290,780){\circle{3}}
\put(290,775){\circle{3}}
\put(290,770){\circle{3}}
\put(290,765){\circle{3}}
\put(290,760){\circle{3}}
\put(290,755){\circle{3}}
\put(290,745){\circle{3}}
\put(290,750){\circle{3}}
\put(290,740){\circle{3}}
\put(290,735){\circle{3}}
\put(290,730){\circle{3}}
\put(290,725){\circle{3}}
\put(290,720){\circle{3}}

\put(291,890){$v_0$}

\put(240,880){\mbox{$\Gamma_1$}{}}
\put(280,850){\mbox{$\Gamma_4$}{}}
\put(380,890){\mbox{$\Gamma_3$}{}}
\put(330,900){\mbox{$\Gamma_2$}{}}
\put(330,870){\mbox{$B=\gamma$}{}}

\put(302,920){\mbox{$\Gamma_5$}{}}

\thicklines
\end{picture}
\end{center}
\caption{ {\em  The right figure shows the events of ${\cal D}_1(\gamma)$ and ${\cal D}_2(\gamma)$
on $L_n(\gamma)$, where the solid circles are open vertices and
the circles are closed vertices.
On the existence of open crossing $\gamma$, we can build five arm paths at $v_0$. $T(\gamma)$ is enclosed
by $L$, $T$, $R$, and $B=\Gamma$.}}
\end{figure}
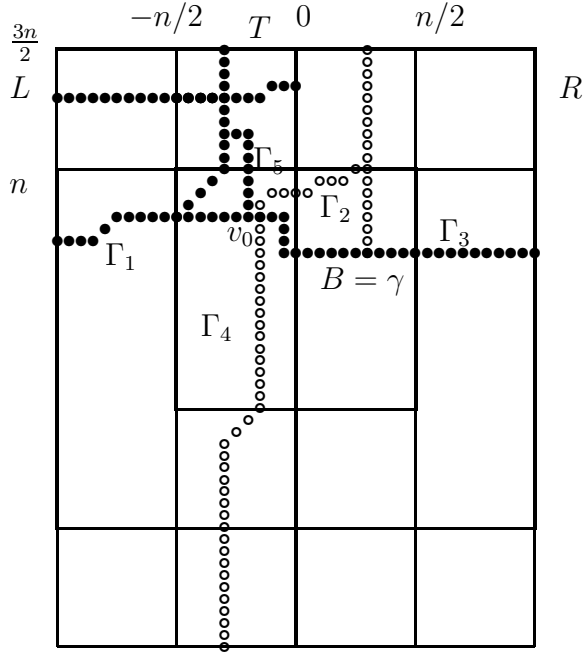

{\bf Proof of Theorem 1.} 
We first estimate the pivotal sites  when $p=0.5$. Let (see Fig 1.)
 Let $L_n(\gamma)$ be the event that the lowest open crossing 
on $[-n, n]^2$ is $\gamma$ for a fixed crossing $\gamma$.
On $L_n(\gamma)$, the lowest crossing $\gamma$ has to stay in $[-n,n]^2$.
Next, for each fixed lowest crossing on $[-n,n]^2$, let (see Fig. 1)
${\cal D}_1(\gamma)$ be the event that there exists two open paths: one from the top of $[-n/2, 0]\times [-n,3n/2]$
to $\gamma$  inside
$T^\circ(\gamma)\cap [-n/2,0]\times [-n, 3n/2]$, and the other one from the left to to the right  in 
$[-n, 0]\times [n, 3n/2]$. In addition, let ${\cal D}_2(\gamma)$ be the event that there exists 
a closed path from the top boundary of $[0,n/2]\times [-n, 3n/2]$ to $\gamma$ inside $T^\circ(\Gamma)\cap [0, n/2]\times [-n, 3n/2]$.  By independent property, for each fixed crossing,
$${\bf P}_p( {\cal D}_1(\gamma)\cap {\cal D}_2(\gamma)\,\,|\,\, L_n(\gamma)) ={\bf P}_p( {\cal D}_1(\gamma)\cap {\cal D}_2(\gamma)).$$
If there exists an open path from 
 the top to the bottom of $[-n/2, 0]\times [-n,3n/2]$ inside of $[-n/2, 0]\times [-n,3n/2]$, then 
there exists an open path from the top of $[-n/2, 0]\times [-n,3n/2]$
to $\gamma$  inside
$T^\circ(\gamma)\cap [-n/2,0]\times [-n, 3n/2]$. With this observation, the FKG ineqaulity, and (2.2),
there exists a constant $C_4>0$ such that
$${\bf P}_{0.5} ({\cal D}_1(\gamma))\geq C_4.$$
The same argument implies that
$${\bf P}_{0.5} ({\cal D}_2(\gamma))\geq C_5$$
for some constant $C_5>0$.
Note that for each fixed $\gamma$, ${\cal D}_1(\gamma)$ and ${\cal D}_2(\gamma)$ are independent, so
there exists $C_6 >0$ such that
$${\bf P}_{0.5} \left({\cal D}_1(\gamma) \cap {\cal D}_2(\gamma)\,\,|\,\, L_n(\gamma)\right)\geq C_6.\eqno{(2.3)}$$
On $L_n(\gamma)$,  the lowest crossing  on $[-n,n]^2$ is $\gamma$.
Note that the boundary of $T(\gamma)$ is a circuit enclosed by  the following four pieces (see Fig. 1):
 the bottom ($B=\gamma$), the top $T$ (the top boundary of $[-n,n]^2$), the left $L$ (the part of the left boundary 
of $[-n ,n]^2$),  and the right $R$ (the part of the right boundary of $[-n, n]^2$). 

With these $L$, $T$, $R$, and $B$, we consider the  left-most top-bottom closed crossing in $T(\gamma)$. 
On $L_n(\gamma)$, if ${\cal D}_1(\gamma) \cap {\cal D}_2(\gamma)$ occurs,
this left-most top-bottom closed crossing $\Gamma_4$ (see Fig. 1)
exists in $T^\circ(\gamma)$ with a starting vertex $v_0\in [-n/2,n/2]\times [-n, n]$ (not included on  $\Gamma_4$) at $\gamma$.
On $L_n(\gamma)$, if  ${\cal D}_1(\gamma) \cap {\cal D}_2(\gamma)$ occurs, we denote by ${\cal E}_{v_0}(\gamma)$
the event that there exists the lowest left-right open crossing $\gamma$ passing through $v_0$, and  there exists 
the  left-most top-bottom closed crossing $\Gamma_4$ from $\gamma$ to the top boundary of $[-n, n]\times [-n, 3n/2]$ with
the starting vertex  $v_0\in [-n/2,n/2]\times [-n,n]$. 
Since $\Gamma_4$ is the  left-most closed crossing, by  the three-arm-argument 
 there exists an additional  open path $\Gamma_5$ inside $T^\circ (\gamma)$
from a neighbor of $v_0$, denoted by  $v_1$,  to the left of $[-n,n]\times [-n, 0, 3n/2]$ (see Fig. 1).
Here $\Gamma_5$ includes $v_1$.

In summary,  if ${\cal E}_{v_0}(\gamma)$ occurs,
there are four disjoint paths from $v_0$ (not including $v_0$):
two  open paths  $\Gamma_1$ and $\Gamma_3$ from $v_0$ to the left and to the right of  $[-n,n]\times [-n, 3n/2]$, 
 and two closed paths $\Gamma_2$ and
$\Gamma_4$ from $v_0$ to the top and to the bottom of $[-n,n]\times [-n, 3n/2]$, respectively.
Also, $v_0$ is open.
In addition to these four arm paths, there exists an open path $\Gamma_5$ in $T^\circ (\gamma)$ from $v_1$
to the left of $[-n,n]\times [-n, 3n/2]$ and $v_1$ is open. Note that $v_1$ is a neighbor of $v_0$.
Note also that if $v_0$ is fixed,  then there are at most nine choices for choosing $v_1$.
Let $N([-n,n]\times [-n, 3n/2])$ be the number of pivotal sites for the crossing $LR([[-n, n]\times [-n, 3n/2])$.
Therefore, for fixed crossing $\gamma$, by  (2.3), Reimer's inequality (2000),  translation invariance,   Lemma 1,
and independent property,
\begin{eqnarray*}
C_6&\leq &{\bf P}_{0.5} \left({\cal D}_1(\gamma) \cap {\cal D}_2(\gamma)\,\,|\,\, L_n(\gamma)\right )\\
&\leq &  {\bf P}_{0.5} \left( \exists \,\, v_0\in [-n/2,n/2]\times [-n, n] \mbox{ such that  ${\cal E}_{v_0}(\gamma)$ occurs},\,\,|\,\, L_n(\gamma)\right)\\
&= &  {\bf P}_{0.5} \left( \bigcup_{v_0\in [-n/2,n/2]\times [-n, n] }{\cal E}_{v_0}(\gamma)\,\,|\,\,L_n(\gamma)\right)\\
&\leq & \sum_{v_0\in [-n/2 ,n/2]\times [-n, n]} {\bf P}_{0.5} \left( {\cal E}_{v_0} (\gamma)\,\,|\,\,L_n(\gamma)\right)\\
&\leq & 9\sum_{v_0 \in [-n/2, n/2]\times [-n, n]} {\bf P}_{0.5} \left( {\bf C}({\bf 0})\cap \partial [-n/2,n/2]^2\neq \emptyset\right) \\
&&\hskip 3cm \cdot {\bf P}_{0.5} \left(\exists \mbox{ four arm paths at $v_0$ for } [-n, n]\times [-n, 3n/2]\,\,|\,\,L_n(\gamma)\right)\\
&\leq & 9 n^{-C_4} {\bf E}_{0.5}\left( N([-n,n]\times [-n, 3n/2, 3n]\,\,|\,\,L_n(\gamma)\right).\hskip 2.2in (2.4)
\end{eqnarray*}
On $L_n(\gamma)$, the lowest open  crossing on $[-n,n]\times [-n, 3n/2]$ stays inside $[-n, n]^2$. Thus, 
on $L_n(\gamma)$, each 
pivotal site for the left-right open crossing of $[-n,n]\times [-n, 3n/2]$ is also a pivotal site for the left-right
open crossing of 
$[-n,n]^2$. In other words, for each $p$, 
$${\bf E}_{p}\left( N([-n,n]\times [-n, 3n/2]\,\,|\,\,L_n(\gamma)\right)\leq {\bf E}_{p}\left( N_n\,\,|\,\,L_n(\gamma)\right).\eqno{(2.5)}$$
Together with (2.4) and (2.5), we have for each fixed crossing $\gamma$,
$$n^{\alpha } \leq {\bf E}_{0.5}\left( N_{n}\,\,\,|\,\,\, L_n(\gamma)\right).\eqno{(2.6)}$$

Now we show Theorem 1 by (2.6).
$${\bf E}_{p}\left( N_{n}; LR([-n, n]^2)\right)=\sum_{\gamma}{\bf E}_{p}\left( N_{n}\,\,\,|\,\,\, L_n(\gamma)\right)
{\bf P}_p ( L_n(\gamma)).\eqno{(2.7)}$$
On $L_n(\gamma)$, if there exist four arm paths at $v$, then by the three-arm-path argument, $v$ is on the lowest crossing.
Therefore, $N_{n}$ is the number of vertices $\{v\}\subset \gamma$ 
such that
there exist closed paths inside $T^\circ(\gamma)$ from $v$ (not including $v$) to the top of $[-n, n]^2$.
For each fixed crossing $\gamma$, let $V_n(\gamma)$ be the vertices of $\{v\}$ above. 
By the independence property of the lowest crossing, these closed paths only depend on the configurations on $T^\circ(\gamma)$: 
$${\bf E}_{p}\left( N_{n}\, \,\,|\,\,\, L_n(\gamma)\right)={\bf E}_{p}\left( V_n(\gamma)\right).$$
Note that ${\bf E}_p V_n(\gamma)$ is  decreasing  in $p$ for each fixed crossing $\gamma$.
Therefore,
$${\bf E}_{0.5}\left( N_{n}\, \,\,|\,\,\, L_n(\gamma)\right)={\bf E}_{0.5}\left( V_n(\gamma)\right)\leq {\bf E}_{p}\left( V_n(\gamma)\right)={\bf E}_{p}\left( N_{n}\,\,\,|\,\,\, L_n(\gamma)\right).\eqno{(2.8)}$$
By (2.6), (2.7), and (2.8), for all $p \leq 0.5$,
$$n^\alpha {\bf P}_p (LR([-n,n]^2))\leq {\bf E}_{p}\left( N_{n}; LR([-n,n]^2)\right).\eqno{(2.9)}$$
Theorem 1 follows. $\Box$\\

{\bf Remark.} Kesten, Sidoravicius, and Zhang (1998) gave a precise order of the 
probability estimate for the five arm paths.
The proof is quite  long. \\

If we denote that $N^*_n$ be the pivotal sites for a closed crossing in $[-n,n]^2$, then by symmetry and Theorem 1, 
we have
the following Corollary.\\

{\bf Corollary 4.} {\em If $q \geq 0.5$, then}
$${\bf E}_{q} \left( N_n^* ; LR^*([-n,n]^2)\right) \geq n^{\alpha}.\eqno{(2.10)}$$

{\bf Proof of Corollary 2.} 
By Theorem 1 and Russo's formula (see (2.30) in Grimmett (1999)), 
note that $LR([-n,n]^2)$ is an increasing  event, so there exist $C_i=C_i(p)$ for $i=7,8$ such that for $p < 0.5$,
$$ {\bf P}_{p} (LR([-n,n]^2)) \leq \exp\left( -\int_{p}^{0.5} {\bf E}_p (N_n\,\,|\,\, LR([-n,n]^2)\right)\leq C_7\exp(- C_8n^{\alpha}).\eqno{(2.11)}$$
By (2.10) and symmetry,
if $q >0.5$, then
$$ {\bf P}_{q} (LR^*([-n,n]^2)) \leq C_7\exp(- C_8n^{\alpha}).\eqno{(2.12)}$$
Note that if $|{\bf C}({\bf 0}) |\geq n$, then there exists an open path from the origin
to $\partial [-\sqrt{n} ,\sqrt{n}]^2$. By (2.12), symmetry and the FKG inequality,  
there exist $C_i=C_i(p)$ for $i=9,10$ such that for $p < 0.5$,
$${\bf P}_{p} (|{\bf C}({\bf 0}) |\geq n) \leq [{\bf P}_p( LR([-\sqrt{n}, \sqrt{n}]^2))]^{1/2} \leq C_9 \exp(- C_{10}n^{\alpha/2}).\eqno{(2.13)}$$
By (2.13),
$$p_T\geq 0.5.\eqno{(2.14)}$$
Corollary 2 also follows from  (2.13) and a simple computation (see Theorem 5.4 in Grimmett (1999)). $\Box$\\

{\bf Proof of Corollary 3.}
By Lemma 1, for each $n$,
$${\bf P}_{0.5} \left( |{\bf C}({\bf 0}) |=\infty\right)\leq {\bf P}_{0.5} \left({\bf C}({\bf 0})\cap \partial [-n,n]^2\neq \emptyset\right )\leq 
n^{-C_4}.\eqno{(2.15)}$$
Thus, $\theta(0.5)=0$, so 
$$p_c \geq 0.5.\eqno{(2.16)}$$
Now we assume that $p_c > 0.5$ and select $0.5 < q < p_c  $.
With this assumption, by (2.2),
$${\bf P}_q ( LR^*([-n,n]^2)) \geq C_3>0.\eqno{(2.17)}$$
 Since (2.12) and (2.17) cannot hold together for large $n$, the contradiction tells us that $p_c \leq 0.5$.
Together with (2.16), we have  $p_c=0.5$.
Note that $p_T\leq p_c$, so by (2.14), $p_T=0.5$. Therefore,
Corollary 3 follows. $\Box$\\

\begin{center}{\large \bf References} \end{center}
Grimmett, G. (1999). {\em Percolation}. { Springer-Verlag},  New York.\\
Harries, E. (1960). A lower bound for the critical probability in a certain percolation process. 
{\em Proceedings of the Cambridge Philosophical Society} {\bf 56} 13--20.\\
Kesten, H. (1980). The critical probability of bond percolation on the square lattice equals 1/2. {\em Comm. Math. Phys.}
{\bf 74} 41--59.\\
Kesten, H. (1982).  {\em Percolation Theory for Mathematicians}, { Birkhauser}, Boston.\\
Kesten, H., Sidoravicius, V. and  Zhang, Y. (1998). Almost all words are seen in critical site percolation on the triangular lattice. {\em Electron. J. Probab.} {\bf 3} 1--75. \\
Reimer, D. (2000). Proof of the van den Berg-Kesten inequality.  {\em Combin. Probab. Compute} {\bf 9} 27--32.\\
Russo, L. (1978). A note on percolation.
\emph{Z. Wahrsch. verw. Gebiete} {\bf 43} 39--48.\\
Seymour, P. D. and Welsh, D. J. A. (1978).
Percolation probabilities on the square lattice. In  {\em Advances in Graph Theory} (B. Bollobas ed.) 227--245.
{Ann. of Discrete Math.}  {\bf 3}, North-Holland, Amsterdam.\\

Yu Zhang\\
Department of mathematics\\
University of Colorado\\
Colorado Springs, CO 80933\\
yzhang3@uccs.edu

\end{document}